\documentclass[11pt,article]{article}

\usepackage{amsmath}
\usepackage{amsthm}
  \usepackage{paralist}
  \usepackage{graphics} %% add this and next lines if pictures should be in esp format
  \usepackage{epsfig} %For pictures: screened artwork should be set up with an 85 or 100 line screen
\usepackage{graphicx}
\usepackage{caption}
\usepackage{subcaption}
\usepackage{epstopdf}%This is to transfer .eps figure to .pdf figure; please compile your paper using PDFLeTex or PDFTeXify.
 \usepackage[colorlinks=true]{hyperref}
 \usepackage{multirow}
\input{amssym.tex}

\newtheorem{theorem}{Theorem}[section]
\newtheorem{corollary}{Corollary}
\newtheorem{lemma}[theorem]{Lemma}
\newtheorem{proposition}{Proposition}

 \numberwithin{equation}{section}

\newcommand{\keywords}

\def\bc{\begin{center}}       \def\ec{\end{center}}
\def\ba{\begin{array}}        \def\ea{\end{array}}
\def\be{\begin{equation}}     \def\ee{\end{equation}}
\def\bea{\begin{eqnarray}}    \def\eea{\end{eqnarray}}
\def\beaa{\begin{eqnarray*}}  \def\eeaa{\end{eqnarray*}}

\def\mathbb{\Bbb}

\begin{document}

\title{\bf Boundary spikes of a Keller--Segel chemotaxis system with saturated logarithmic sensitivity}
\author{Qi Wang \thanks{(Email:{\tt qwang@swufe.edu.cn}).  }\\
Department of Mathematics\\
Southwestern University of Finance and Economics\\
555 Liutai Ave, Wenjiang, Chengdu, Sichuan 611130, China
}

\date{}
\maketitle

%
%\title{\LARGE \bf On a nonlinear boundary value problem with integral constraint arising from species competitions}

%
%\date{}
%\maketitle

%\title[Nonlinear boundary value problem with integral constraint]{\LARGE \bf On a nonlinear boundary value problem with integral constraint arising from species competitions}
%
%
%\centerline{\scshape   }
%
%{\footnotesize
%% please put the address of the first author
% \centerline{Department of Mathematics }
%   \centerline{Southwestern University of Finance and Economics}
%   \centerline{555 Liutai Ave, Wenjiang, Chengdu, Sichuan 611130, China}
%} % Do not forget to end the {\footnotesize by the sign }

%
%\title{\LARGE \bf Qualitative analysis of stationary Keller-Segel chemotaxis models with logistic growth}
%\author{Qi Wang \thanks{Department of Mathematics, Southwestern University of Finance and Economics, Chengdu, Sichuan 611130, China
%        ({\tt qwang@swufe.edu.cn}).}}
%\date{}

\abstract
In this paper, we study the nonconstant positive steady states of a Keller-Segel chemotaxis system over a bounded domain $\Omega\subset \mathbb{R}^N$, $N\geq 1$.  The sensitivity function is chosen to be $\phi(v)=\ln (v+c)$ where $c$ is a positive constant.  For the chemical diffusion rate being small, we construct positive solutions with a boundary spike supported on a platform.  Moreover, this spike approaches the most curved part of the boundary of the domain as the chemical diffusion rate shrinks to zero.  We also conduct extensive numerical simulations to illustrate the formation of stable boundary and interior spikes of the system.  These spiky solutions can be used to model the self--organized cell aggregation phenomenon in chemotaxis.

\textbf{Competition model, shadow system, nonlinear boundary value problem, transition layer.}

1
\section{Introduction and main result}

Chemotaxis is the oriented movement of cells along the gradient of certain chemicals in their environment.  One of the most interesting phenomena in chemotaxis is the cell aggregation and it is an important mechanism for the formation of a fruiting body from single cells.  See \cite{DW, Ho} and \cite{Ho2} for detailed discussions.  This phenomenon has been discovered in bacteria, in particular E. coli, and also in the cell slime molds such as Dictyostelium discoideum.  See \cite{BWS, DW} and the references therein.  One of the main reasons of the directed movements of the cells is that they can navigate within a complex environment by detecting and processing of a variety of internal and external signals.  Certain chemicals, most frequently inorganic salt, amino acids and some proteins called chemokines, can induce chemotaxis in motile cells.

In this paper we study the existence of nonconstant positive solutions $(u,v)=(u(x),v(x))$ to the following system
\begin{equation}\label{6}
\left\{
\begin{array}{ll}
\nabla\cdot (  \nabla u - p u \nabla \ln (v+c))=0, & x \in\Omega,\\
\epsilon^2 \Delta v - v+\beta u/\alpha=0, & x \in\Omega,\\
\frac{\partial u}{\partial \textbf{n}}=\frac{\partial v}{\partial \textbf{n}}=0, & x \in \partial \Omega,
\end{array}
\right.
\end{equation}
where $\Omega$ is a bounded domain in $\mathbb R^N$, $N\geq1$, with smooth boundary $\partial \Omega$.  $\textbf{n}$ is the unit outer normal on the boundary.  The parameters $p$, $c$, $\epsilon$, $\alpha$ and $\beta$ are assumed to be positive constants.  Let $u$ be the cellular population density and $v$ be the chemical concentration, then system (\ref{6}) can be used to model the steady state of a chemotaxis system with saturated logarithmic sensitivity function.  We are concerned with boundary spike solutions to (\ref{6}) and our main theorem goes as follows.
\begin{theorem}\label{theorem11}
Suppose that $p\in(1,\infty)$ for $N=1,2$ and $p\in(1,(N+2)/(N-2))$ for $N\geq 3$.  Assume that $\int_\Omega u(x)dx=M\leq \frac{\alpha c \vert \Omega \vert}{\beta (p-1)}$, where $\vert \Omega \vert$ is the $N$--dimensional volume of domain $\Omega$.  Then for each $\epsilon>0$ small, there exists a positive nontrivial solution $(u_\epsilon(x),v_\epsilon(x))$ to (\ref{6}) such that
\[u_\epsilon(x)=\frac{M(v_\epsilon(x)+c)^p }{\int_\Omega (v_\epsilon(x)+c)^p dx}.\]
For $\epsilon$ being sufficiently small, there exists a small region $\Omega_\epsilon^{(i)} \subset \Omega$ such that $v_\epsilon(x)$ achieves its global maximum at a boundary point $P_\epsilon \in \partial \Omega_\epsilon^{(i)} \cap \partial \Omega$.  Moreover, $v_\epsilon(x)$ has a boundary spike supported on a platform with height $t_\epsilon^*>0$ in the following sense:

(i) $\text {diam} (\Omega_\epsilon^{(i)}) \leq C \epsilon$, where $C$ is a positive constant independent of $\epsilon$;

(ii) $\max_{\bar{\Omega}}v_\epsilon (x)=v_\epsilon (P_\epsilon)>\frac{\beta M}{\alpha \vert \Omega \vert}$;

(iii)  the profile of $v_\epsilon$ can be approximated as
\[v_\epsilon(x)=\mathbf{w}_\epsilon\left( \frac{\vert x-P_\epsilon \vert}{\epsilon} \right)+t^*_\epsilon+o(\epsilon),x\in \Omega,\]
where
\begin{equation}\label{7}
t^*_\epsilon\in \Big(0,\frac{c}{p-1}\Big)\text{ and }t^*_\epsilon \rightarrow \frac{\beta M}{\alpha \vert \Omega \vert} \text{ as } \epsilon \rightarrow 0,
\end{equation}
and $\mathbf{w}_\epsilon$ is a radially symmetric function satisfying

(iii a) $\mathbf{w}_\epsilon(x)>0$, $x\in \mathbb{R}^N$ and $\mathbf{w}_\epsilon(0)=\max_{\mathbb{R}^N}\mathbf{w}_\epsilon(x)$,

(iii b) $\mathbf{w}'_\epsilon(r)<0$, $r=\vert x \vert$, for $x \in \mathbb{R}^N\backslash\{0\}$,

(iii c) $\mathbf{w}_\epsilon(r) \leq Ce^{-\mu r}$, $r>1$, for some positive constants $C$ and $\mu$ independent of $\epsilon$.
Furthermore, $P_\epsilon$ approaches the most curved part of the boundary $\partial \Omega$ in the sense that $H(P_\epsilon)\rightarrow \max_{P\in \partial \Omega}H(P)$,
where $H(P)$ is the mean curvature of the boundary $\partial \Omega$ at point $P$.
\end{theorem}

The mathematical modeling of chemotaxis was initiated by Keller and Segel in their pioneering works \cite{KS,KS1,KS2} during 1970s.  Denote by $u(x,t)$ the cell population density and by $v(x,t)$ the chemical concentration at space--time $(x,t)$.  Then the general Keller--Segel model reads as follows
\begin{equation}\label{1}
\left\{
\begin{array}{ll}
u_t=\nabla\cdot (d_1(u,v) \nabla u - \chi (u,v) \nabla \phi(v)), & x \in\Omega,~t>0,\\
v_t=d_2 (u,v) \Delta v+k(u,v),& x \in\Omega,~t>0,\\
u(x,0)=u_0(x)\geq 0,v(x,0)=v_0(x)\geq 0,& x\in\Omega,\\
\frac{\partial u}{\partial \textbf{n}}=\frac{\partial v}{\partial \textbf{n}}=0, & x \in \partial \Omega,~t>0,\\
\end{array}
\right.
\end{equation}
where $d_1, d_2>0$ are the diffusion rates of cells and chemicals respectively and $\chi(u,v)$ is a function that measures the chemotactic response of cells to the
chemical.  Moreover, $\chi$ is positive if the chemical is chemo--attractive (sugar, nutrition, etc.) to the cells and it is negative if the chemical is chemo--repulsive (poison, toxin for example).  $\phi(v)>0$ is the so--called (signal--dependent) sensitivity function and it reflects the variance of cellular sensitivity with respect to the chemical concentration.  $k(u,v)$ is the creation and degradation rate of the chemical.  The non--flux boundary conditions mean that the domain $\Omega$ is enclosed thus it inhibits cell immigration and chemical diffusion across the boundary.  The nonnegative initial data $u_0,v_0$ are not identically zero.

Keller--Segel system (\ref{1}) and its variants have been studied extensively over the past few decades.  These systems are very rich in mathematical dynamics and they can induce many interesting and striking properties such as global solutions, finite time blow--ups, traveling waves, etc.  Furthermore, the intuitively simple system (\ref{1}), even in its simplest form, successfully demonstrates its ability in presenting solutions with self--organized spatial patterns with concentration properties such as boundary spike or interior spike.  These spiky solutions can be used to model the aforementioned cellular aggregation and concentration phenomena.  See the survey paper \cite{HP} for detailed discussions.

From the viewpoint of mathematical modeling, we say that \emph{cells aggregate} if there exists a solution $u(x,t)$ that converges to a nonconstant positive function $u(x)$,
\[\lim_{t\rightarrow T } u(x,t)=u(x),T\leq \infty.\]
In this paper, we focus on the interplay between the sensitivity function $\phi(v)$ and the terms $d_1$, $d_2$ and $\chi$ on the mathematical modeling of cell aggregations phenomenon.  To make our goal evident, we assume that $d_1$ and $d_2$ are positive constants and $k(u,v)=-\alpha v +\beta u$ for some positive constants $\alpha$ and $\beta$.  Under these settings, Keller--Segel model (\ref{1}) becomes
\begin{equation}\label{2}
\left\{
\begin{array}{ll}
u_t=\nabla\cdot (d_1 \nabla u - \chi u \nabla \phi(v)), & x \in\Omega, ~t>0,\\
v_t= d_2 \Delta v - \alpha v+\beta u, & x \in\Omega,~ t>0,\\
u(x,0)=u_0(x)\geq 0,~v(x,0)=v_0(x)\geq 0,& x\in \Omega,\\
\frac{\partial u}{\partial \textbf{n}}=\frac{\partial v}{\partial \textbf{n}}=0, & x \in \partial \Omega,~t>0,
\end{array}
\right.
\end{equation}
and we want to find nonconstant positive steady states to systems in the form of (\ref{2}).

Time--dependent system (\ref{2}) can model cell aggregation if $L_\infty$--norm of the solutions goes to infinity, then the aggregation is simulated by a $\delta$--function or a linear combination of $\delta$-functions that represents the chemical concentration.  The formal study of blow--up solutions of Keller-Segel models was initiated in \cite{CP,V} and the local profile of the blow--ups was rigourously constructed in \cite{HV} over a 2D domain.  Though such blow--up solutions are evidently connected to the phenomenon of cellular aggregation, a $\delta$--function is not the optimal choice from the viewpoint of mathematical modeling.  It challenges the rationality that chemical concentration can not be infinity at a spot in the domain and it also brings difficulties to numerical simulations.  Moreover, rigourous analysis of the states after blow--up is impossible.  An alternative approach available in the literature is to show that positive solutions of (\ref{2}) exist globally in time and converge to bounded steady states as time increases.  Then the steady states with concentrating properties such as spikes, transition layers, etc. can be used to describe the cell aggregation phenomenon.  We shall take the second approach in our paper and it is the goal of the current work to investigate the existence of nonconstant positive steady states of system (\ref{2}).

For $N=1$, it has been well known that solutions to (\ref{2}) and a collection of its variants are global and bounded in time.  See \cite{HPS, OY} for example.  However, the global existence is not always available for higher dimensional domains.  It turns out that the sensitivity function $\phi(v)$ and system parameters in (\ref{2}) play essential roles in existence of its global--in--time solutions.  Two most commonly used sensitivity functions are $\phi(v)=v$ and $\phi=\ln v$, which lead (\ref{2}) to the so--called minimal model and logarithmic model respectively.  See the survey paper of Hillen and Painter \cite{HP}.

For $\phi(v)=v$, Nagai, \emph{et al}. \cite{NSY} established the global solutions $(u(x,t),v(x,t))$ of (\ref{2}) with $d_1=d_2=1$ over $\Omega \subset \mathbb{R}^2$, if either $\int_\Omega u_0(x)dx < 4\pi/(\beta \chi)$, or $\int_\Omega u_0(x)dx < 8\pi/(\beta \chi)$, $\Omega$ is a disk and $(u_0(x),v_0(x))$ is radial in $x$.  Winkler \cite{Wk} studied the same problem with $d_1=d_2=\chi=\alpha=\beta=1$ over $\Omega \subset \mathbb{R}^N,N\geq3$.  It has been proved that the solutions of (\ref{2}) exist globally in time provided that $\Vert u_0(x) \Vert_{L^{N/2+\delta}}$ and $\Vert \nabla v_0(x) \Vert_{L^{N+\delta}} $ are small for some $\delta>0$.

For $\phi(v)=\ln v$, Nagai, \emph{et al}. \cite{NSY2} established the global solutions of (\ref{2}) with $d_1=d_2=\alpha=\beta=1$ over $\Omega \subset \mathbb{R}^2$ provided that either $\chi<1$ or $\chi<5/2$, $\Omega$ is a disk and $(u_0(x),v_0(x))$ is radial in $x$.  Winkler \cite{Wk1} studied the same problem and established the global classical solutions over $\Omega \subset \mathbb{R}^N,N\geq 2$ if $\chi<\sqrt{2/N}$ and global weak solutions if $\chi<\sqrt{(N+2)/(3N-4)}$.  Furthermore, Nagai and Senba \cite{NS} proved the global existence and boundedness for the parabolic--elliptic system of (\ref{2}) when $\Omega \subset \mathbb{R}^2$ is a disk and $(u_0(x),v_0(x))$ is radial, or $\chi<2$, $\Omega \subset \mathbb{R}^3$ is a ball and $(u_0(x),v_0(x))$ is radial.  An exploration of literatures suggests that both the smallness of initial data and of the chemo--attraction rate $\chi$ tend to prevent solutions of (\ref{2}) from blowing up in finite or infinite time.

Attempts have been made to show that the steady states of (\ref{2}) have concentrating properties such as boundary spikes, interior spikes, etc., when the chemical diffusion rate is sufficiently small.  This approach was initiated by Lin, Ni and Takagi in a series of remarkable papers \cite{LNT, NT, NT2} around 1990s.  To be precise, their pioneering works focus on the logarithmic case $\phi(v) =\ln v$ and can be summarized as follows.  Through the transformation introduced in \cite{Sc} and proper parameter scalings, the stationary system of (\ref{2}) reduces to the following single equation
\begin{equation}\label{3}
d_2 \Delta v -v +v^p=0, x\in \Omega,~~\nabla v \cdot \textbf{n} =0,~x\in \partial \Omega,
\end{equation}
where $p=\frac{\chi}{d_1}$.  Equation (\ref{3}) is associated with an energy functional
\[J_{d_2}:=\frac{1}{2}\int_\Omega (d_2 \vert \nabla v \vert^2+v^2)dx-\frac{1}{p+1} \int_\Omega v^{p+1}_{+} dx\]
with $v^+=\max\{v,0\}$.  Suppose that $p\in(1,\infty)$ if $N=1,2$ and $p\in(1,\frac{N+2}{N-2})$ if $N\geq 3$, Lin, Ni and Takagi \cite{LNT} applied variational techniques and min--max arguments to prove that if $d_2$ is small, (\ref{3}) has only nonconstant positive least--energy solution, which is a critical point of $J_{d_2}$ in $H^1(\Omega)$. Furthermore, Ni and Takagi \cite{NT,NT2} showed that the least--energy solutions must concentrate at a boundary spike if $d_2$ is sufficiently small and the spike moves to the place of $\partial \Omega$ with the largest mean curvature as $d_2$ goes to zero.   Motivated by their works, Grossi, Pistoia and Wei \cite{GPW} obtained solutions of (\ref{3}) with a single interior spike that stays in the most centered region of the domain (in the sense of a critical point of the distance function) for small $d_2$.  By the Lyapunov-Schmidt method, Gui and Wei \cite{GW} construct multi-peak solutions with $m$ interior spikes and $n$ boundary spikes, for any nonnegative integers $m$, $n$, as long as $d_2$ is sufficiently small.  The methods and techniques used in these brilliant works began with good guesses on the structures of the spike(s).  Moreover, the stability analysis has also been carried out through the infinite dimension Lyapunov--Schmidt reduction and the method of NLEP (nonlocal eigenvalue problem) recently developed by J. Wei, \emph{et al}.  See \cite{SWW} and the surveys in \cite{N,N2} for recent development in this direction.

It is worthwhile to point out that a large chemotactic coefficient $\chi$ can also destabilize the constant solution of (\ref{2}) like a small diffusion rate.  Therefore, nontrivial steady states may emerge through large chemotactic coefficient $\chi$.  For example, by taking $\chi$ as the bifurcation parameter, Wang, \emph{et al}. \cite{CKWW,W,WX} applied global bifurcation methods to investigate nontrivial steady states of (\ref{2}) over one--dimensional finite interval.  Furthermore, these solutions are shown to admit boundary spikes or interior transition layers as $\chi$ approaches infinity.

The logarithmic sensitivity function $\phi(v)=\ln v$ was chosen largely due to the Weber--Fechner's law describing cellular behaviors: \emph{the subjective sensation is proportional to the logarithm of the stimulus intensity} \cite{HP}.  For the logarithmic model (\ref{2}) with $\phi(v)=\ln v$, we observe that the dynamics of cellular movements are dominated by the taxis flux $\frac{\chi}{v} \nabla v$ which may be unbounded for $v \approx 0$.  We want to note that this singularity is an important mechanism in the formation of traveling wave solutions as showed in \cite{KS} or discussed in the survey paper \cite{ZW}.  See \cite{LW,TWW} for more works on the logarithmic model.  However, from the viewpoint of mathematical modeling, it is not reasonable to assume that low chemical concentration elicits significant chemotactic responses from the cell's motility, hence we modify this problem by taking $\phi(v)=\ln(v+c)$ for a positive constant $c$ which saturates $\phi'(v)$ at $v=0$.  Then the steady state of (\ref{2}) becomes
\begin{equation}\label{4}
\left\{
\begin{array}{ll}
\nabla\cdot (d_1 \nabla u - \chi u \nabla \ln (v+c))=0, & x \in\Omega,\\
d_2 \Delta v - \alpha v+\beta u=0, & x \in\Omega,\\
\frac{\partial u}{\partial \textbf{n}}=\frac{\partial v}{\partial \textbf{n}}=0, & x \in \partial \Omega,
\end{array}
\right.
\end{equation}
where $d_1$, $d_2$, $\chi$, $c$, $\alpha$ and $\beta$ are positive constants.  First of all, we claim that the total cell population in the time--dependent counterpart of (\ref{4}) is conserved.  Indeed, integrating the first equation over $\Omega$, we have from the Neumann boundary conditions that
\[\frac{d}{dt} \int_\Omega u(x,t)dx=\int_\Omega u_t(x,t)dx=\int_\Omega \nabla\cdot \Big(d_1 \nabla u - \chi u \nabla \ln (v+c)\Big)dx =0,\]
and this implies that for all $t>0$,
\begin{equation} \label{5}
\int_\Omega u(x,t) dx=\int_\Omega u(x,0) dx=M.
\end{equation}
Throughout this paper, we assume the cell population in (\ref{5}) is a fixed positive constant denoted by $\int_\Omega u(x) dx =M$.  Moreover, since small chemical diffusion rate supports nontrivial positive solutions to (\ref{4}) as we shall show, we denote for the simplicity of notations that
\[\epsilon^2=\frac{d_2}{\alpha}, \text{~and~} p=\frac{\chi}{d_1},\]
then it is easy to see that (\ref{4}) leads us to (\ref{6}).

According to our results in Theorem \ref{theorem11}, $v_\epsilon(x)$ concentrates at the boundary point $P_\epsilon$ if the chemical diffusion rate $\epsilon$ is sufficiently small.  Compared to the spiky solutions in \cite{LNT,N,N2} and \cite{SWW}, $v_\epsilon(x)$ has a boundary spike which is supported on a platform with height $t^*_\epsilon$ in (\ref{7}).  The results in \cite{LNT,N,N2} and \cite{SWW} are based on the application of the classical variational methods due to Ambrosetti and Rabinowitz \cite{AR} on the following equation,
\[d_2 \Delta v -v +f(v)=0, x\in \Omega; \nabla v \cdot \textbf{n} =0,~x\in \partial \Omega,\]
equipped with the energy functional
\[J_{d_2}:=\frac{1}{2}\int_\Omega (d_2 \vert \nabla v \vert^2+v^2)dx- \int_\Omega F(v)dx,~F(t)=\int_0^t f(s)ds.\]
It is assumed, among other conditions, that $f(t)\equiv 0$ for $t\leq0$ and $f(t)=o(t)$ as $t\rightarrow 0^+$.  As we shall see later, system (\ref{6}) can be reduced to the single equation (\ref{10}), however, $f(v)$ there apparently does not satisfy these conditions, therefore it inhibits the direct application of the results in \cite{LNT,N,N2}.  On the other hand, one needs some nontrivial modifications of the variational methods in \cite{AR} in order to apply them for (\ref{10}).  Furthermore, the nonlocal structure in (\ref{10}) can not be scaled out as in \cite{LNT} and this makes the application of previous results on (\ref{10}) difficult.  Therefore the adaption of the ideas in \cite{LNT}, etc. to (\ref{6}) requires several technical tricks.  Compared to the results obtained in \cite{NT,NT2} and \cite{SWW}, where the sensitivity function is $\phi(v)=\ln v$, our boundary spike is supported on a platform with height $t^*>0$.  If $c=0$, we see from (\ref{7}) that our results coincide with those obtained in \cite{NT,NT2}, except that the assumption on the initial cell population becomes $M\leq0$, which apparently can not be true.  However, we can consider an approximation problem with $c\approx0$, then we can choose the parameters properly, for example $\alpha$ being large, such that both the assumption on the initial cell population $M$ and the smallness of $\epsilon$ can be achieved.

We want to point out that the assumptions on $M$ and $p$ are made out of technical reasons.  As we have discussed earlier, either small initial total cell mass or small chemo--attraction coefficient $\chi=pd_1$ tends to prevent solutions of (\ref{2}) from blowing up in finite or infinite time, therefore it makes our analysis of pattern formations in the steady state realistic.  It is also worthwhile to mention that the global existence of a parabolic--elliptic type system (\ref{2}) with $\phi(v)=\ln (v+c)$ was studied in \cite{B}, and the fully parabolic system has been investigated in \cite{QW} for $d_1=d_2=\alpha$ provided that $\chi/d_1=p<\sqrt{2/N}$, $N\geq 2$.  However, the full understanding about its global existence or finite time blow--ups is away from complete and there are many problems widely left open.  See the surveys \cite{HP,Ho,Ho2}.  Though our assumption on $p=\chi/d_1$ in Theorem \ref{theorem11} does not guarantee the global existence of the time--dependent systems, we surmise that the smallness assumption of the cell population might contribute in this matter.

The rest of this paper is organized as follows.  In Section \ref{section2}, we convert system (\ref{6}) into a single equation with a nonlocal structure.  Then we establish the existence of nontrivial positive solutions of the corresponding local problem.  Section \ref{section3} is denoted to study the asymptotic behaviors of the positive solutions obtained in Section \ref{section2} as $\epsilon \rightarrow 0$. The proof of Theorem \ref{theorem11} is included at the end of this section.  Finally, we present numerical simulations in Section \ref{section4} to illustrate and verify our numerical findings.  Throughout the rest of this paper, we assume that $C$ and $C_i$, $i=1,2,...$ are generic positive constants that may vary from line to line.

\section{Existence of nonconstant positive solutions}\label{section2}

First of all, we observe that system (\ref{6}) can be reduced to a single equation with an integral constraint.  Indeed, we see that $u$-equation in (\ref{6}) is equivalent to
 \[\nabla \cdot \left(  u \nabla \Big(\ln u-p \ln (v+c) \Big)\right)=0.\]
Testing it by $\ln u -p \ln (v+c)$ over $\Omega$ by parts, we have that
 \[ \int_\Omega u \Big\vert \nabla \Big(\ln u-p \ln (v+c) \Big)\Big \vert^2dx=0,\]
therefore, there exists a positive constant $C$ such that
\begin{equation} \label{8}
u=C (v+c)^p
\end{equation}
Integrating (\ref{8}) over $\Omega$ leads us to
\begin{equation}\label{9}
u=\frac{M (v+c)^p}{\int_\Omega (v+c)^pdx},
\end{equation}
where $M=\int_\Omega u(x) dx$ is a fixed positive constant representing the total cell population.  Denoting $m=\beta M/\alpha$, we readily see that $v(x)$ satisfies the following nonlocal equation,
\begin{equation}\label{10}
\left\{
\begin{array}{ll}
\epsilon^2 \Delta v-v+\frac{m  (v+c)^p}{\int _\Omega (v+c)^p dx}=0,&x\in \Omega,\\
\frac{\partial v}{\partial \textbf{n}}=0,&x \in \partial \Omega,
\end{array}
\right.
\end{equation}
which serves as a prototype of (\ref{6}).  Indeed, if $v(x)$ is a positive solution to (\ref{10}), then it solves (\ref{6}) together with $u(x)$ obtained through (\ref{9}).

Two main difficulties emerge in proving the existence of positive nonconstant solutions of (\ref{10}).  First of all, the nonlocal structure $\int _\Omega (v+c)^p dx$ inhibits the direct applications of some standard arguments such as degree--theory etc.  If $c=0$, the equation (\ref{10}) can be transformed into the following problem without a nonlocal term
\[\epsilon^2 \Delta w-w+w^p=0,\]
through the scaling $v=mw(\int_\Omega w^p)^{-1}$ as in \cite{NT}, however, such linear scaling is not available to eliminate the nonlocal structure in (\ref{10}).  Moreover, it is well known that (\ref{10}) is naturally associated with an energy functional and this motivates one to apply variational methods to obtain nontrivial solutions, e.g, the Mountain Pass Theorem, \cite{AR,DN}, etc.  However, the presence of $c$ in (\ref{10}) requires nontrivial modification of these variational methods to this end.

Motivated by the ideas initiated by Ni and Takagi in \cite{NT, NT2} and developed by Sleeman \emph{et al.} in \cite{SWW}, we first investigate the solutions of (\ref{10}) with $\int_\Omega (v+c)^p dx$ replaced by a positive constant $\delta$, while we have to treat the dependence of $v(x)$ on $\delta$, then we proceed to find a solution $v(x)$ such that the integral condition $\int_\Omega (v+c)^p dx =\delta$ is satisfied.  To be precise, we first consider the following problem
\begin{equation}\label{11}
\left\{
\begin{array}{ll}
\epsilon^2 \Delta v-v+\frac{m (v+c)^p}{\delta}=0, &x\in \Omega,\\
\frac{\partial v}{\partial \textbf{n}}=0,& x \in \partial \Omega,
\end{array}
\right.
\end{equation}
with $\delta>0$ being an arbitrary but fixed constant.  The following results indicate that $\delta$ can not be arbitrarily small if (\ref{11}) has a positive solution.

\begin{proposition}\label{prop331}
If (\ref{11}) has a positive solution $v(x)$, $\delta$ must have a lower bound such that
\begin{equation}\label{12}
\delta\geq\delta_0=mp \left(\frac{pc}{p-1}\right)^{p-1}.
\end{equation}
Moreover, for $\delta>\delta_0$, the function $R_\delta(t)=-t+\frac{m(t+c)^p}{\delta}$ has exactly two positive roots $t_1(\delta)$ and $t_2(\delta)$ with $0<t_1(\delta)<t_2(\delta)$ such that $t_1(\delta)$ is monotone decreasing and $t_2(\delta)$ is monotone increasing as functions of $\delta$.  Furthermore,
\[t_2(\delta) \leq \left(\frac{\delta}{m}\right)^\frac{1}{p-1}, \forall \delta\geq \delta_0\]
and
\begin{equation}\label{13}
t_1(\delta) \rightarrow 0,~t_2(\delta) \rightarrow \infty, \text{ as } \delta \rightarrow \infty.
\end{equation}
\end{proposition}

\begin{proof} Let $v(x)$ be a positive solution of (\ref{11}), then we integrate the first equation over $\Omega$ to obtain that
\[\int_\Omega -v+\frac{m (v+c)^p}{\delta} dx=  \int_\Omega R_\delta(v) dx=0. \]
Since $R_\delta(0)>0$, we must have that $R_\delta(t)$ is negative at its minimum point over $(0, \infty)$.  In particular, we see that
\[R'_\delta(t)=-1+\frac{mp(t+c)^{p-1}}{\delta},\]
and $R_\delta(t)$ achieves its minimum at the critical point
\[t^*(\delta)=\left(\frac{\delta}{mp} \right)^\frac{1}{p-1}-c.\]
Therefore we have that $R_\delta(t^*)<0$ or equivalently $-t^*+\frac{m(t^*+c)^p}{\delta}<0$.  On the other hand, we observe that $(t^*+c)^{p-1}=\frac{\delta}{mp}$, then the inequality above implies that
\[t^*>\frac{c}{p-1}.\]
After substituting the formula of $t^*$ into $t^*>\frac{c}{p-1}$, we readily have the desired lower bound $\delta_0$ in (\ref{12}).
To establish the upper estimate on $t_2$, we see that
\[t_2=\frac{m}{\delta}(t_2+c)^p >\frac{mt_2^p}{\delta},\]
then $t_2<\left(\frac{\delta}{m}\right)^\frac{1}{p-1}$ follows from the inequality above.  To show the monotonicity of $t_1(\delta)$ and $t_2(\delta)$ in $\delta$, we note that $0<t_1(\delta)<t^*(\delta)<t_2(\delta)$.  Differentiating $R_\delta(t)=0$ with respect to $\delta$ gives rise to
\[\frac{dt_i(\delta)}{d\delta}=\frac{t}{mp(t_i(\delta)+c)^{p-1}-\delta },~i=1,2,\]
then this leads us to the desired monotonicity instantaneously thanks to the fact that $t_1(\delta)<t^*<t_2(\delta)$.  For all $\delta\geq\delta_0$, we see that $t_1(\delta)\leq t_1(\delta_0)$ and $t_1(\delta)$ is uniformly bounded in $\delta$.  Finally, we have
\[t_1(\delta)=\frac{m(t_1+c)^p}{\delta} \rightarrow 0 \text{ as }\delta \rightarrow \infty\]
and
\[t_2(\delta)>t^*(\delta)=\left(\frac{\delta}{mp} \right)^\frac{1}{p-1}-c \rightarrow \infty \text{ as }\delta \rightarrow \infty.\]
The proof of this proposition completes.
\end{proof}

We are interested in positive solutions to (\ref{11}), hence we shall assume $\delta \geq \delta_0$ from now on.  In particular, (\ref{11}) allows only constant positive solution $\bar{v} \equiv \frac{c}{p-1}$ if $\delta =\delta_0$.  For the sake of simplicity, we skip the index $\delta$ in $t_1$ and $t_2$ and denote them as the first and the second root of $R_\delta(t)$ correspondingly, unless otherwise noticed.

To study positive solutions to (\ref{11}), we introduce the transformation
\begin{equation}\label{14}
v(x)=\delta^{\frac{1}{p-1}}w(x)+t_1,
\end{equation}
where $t_1$ is obtained in Proposition \ref{prop331}, then (\ref{11}) becomes
\begin{equation}\label{15}
\left\{
\begin{array}{ll}
\epsilon ^2 \Delta w- c_\delta w+ f_\delta(w) =0,&x\in\Omega, \\
\frac{\partial w}{\partial \textbf{n}}=0, &x\in \partial\Omega,
\end{array}
\right.
\end{equation}
where \[c_\delta=1-\frac{mp}{\delta}\big(t_1+c\big)^{p-1} \in(0,1),\]
and the nonlinear term reads
\begin{equation}\label{21}
f_\delta(w)= m\left( \Big(w+\frac{t_1+c}{\delta^\frac{1}{p-1}}\Big)^p- p\Big(\frac{t_1+c}{\delta^\frac{1}{p-1}}  \Big)^{p-1}w -
\Big(\frac{t_1+c}{\delta^\frac{1}{p-1}}\Big)^p  \right).
\end{equation}
To obtain nonconstant positive solutions $w_{\epsilon,\delta}(x)$ of (\ref{15}), we want to take the variational approach, i.e., to find nontrivial critical points of an energy functional associated with (\ref{15}) in certain Hilbert space.  To this end, we endow Sobolev space $H^1(\Omega)$ with norm
\[\Vert u \Vert_\epsilon=\Big(\int_\Omega \epsilon^2 \lvert \nabla u\rvert^2+ c_\delta u^2 dx\Big)^{\frac{1}{2}},\]
and assign (\ref{15}) the energy functional
\begin{equation}\label{17}
J_{\epsilon,\delta}(w)=\frac{1}{2}\int_\Omega \epsilon^2 \lvert \nabla w\rvert^2+c_\delta w^2 dx -\int_\Omega F_\delta(w) dx, w \in H^1(\Omega),
\end{equation}
where
\[F_\delta(t)=\int_0^t f^+_\delta (s) ds\]
with
\begin{equation*}
f^+_\delta(s)=
\left\{
\begin{array}{ll}
f_\delta (s), & s\geq 0, \\
0, & s<0.
\end{array}
\right.
\end{equation*}
It is well known that any critical point of $J_{\epsilon,\delta}$ is a weak solution of (\ref{15}) in $H^1(\Omega)$; moreover, since $p$ is subcritical, the weak solution is classical by the standard elliptic regularity theories involving the bootstrap argument.

To obtain nonconstant critical points of $J_{\epsilon,\delta}$, we shall apply the well--known Mountain Pass Theorem on the energy functional $J_{\epsilon,\delta}$ (\ref{17}) due to Ambrosetti and Rabinowitz in \cite{AR} and Ni, \emph{et al}. in \cite{DN,NT} which states that,
\[\mathcal{C}_{\epsilon,\delta}:=\inf_{g\in\Gamma} \sup_{t\in [0,1]} J_{\epsilon,\delta}(g(t))=\inf_{v\in H^1(\Omega),v>0} \sup_{t\geq0} J_{\epsilon,\delta}(tv) \]
is the least among all critical values of $J_{\epsilon,\delta}$ in $H^1(\Omega)$, where
\[\Gamma=\{g\in C([0,1];H^1(\Omega)) ~\vert~ g(0)=0,g(1)=e(x) \}\]
and $e(x) \not\equiv 0$ is non--negative with $J_{\epsilon,\delta}(e)=0$.  Moreover, it is known that $\mathcal{C}_{\epsilon,\delta}$ is independent of the choice of $e$.  The critical point $w_{\epsilon,\delta}(x)$ of (\ref{17}) in $H^1(\Omega)$ is called a least--energy solution of (\ref{15}) and $J_{\epsilon,\delta}(w_{\epsilon,\delta})=\mathcal{C}_{\epsilon,\delta}$ is called \emph{the least--energy value}.

To apply the Mountain Pass Theorem, we first observe that $J_{\epsilon,\delta} \in C^1(H^1(\Omega);\mathbb{R})$, $J_{\epsilon,\delta}(0)=0$ and $J_{\epsilon,\delta}$ satisfies the Palais--Smale condition thanks to (\ref{16}) and the proof of Lemma 3.6 in \cite{AR}; moreover, the argument that leads to Lemma 3.1 in \cite{AR} together with (\ref{16}) yields that, there exist some positive constants $r$ and $C$ such that $J_{\epsilon,\delta}(w_{\epsilon,\delta})>0$ if $0<\Vert w_{\epsilon,\delta}\Vert_\epsilon<r $ and $J_{\epsilon,\delta}(w_{\epsilon,\delta})\geq C>0$ if $\Vert w_{\epsilon,\delta}\Vert_\epsilon=r$.  We now need to show, for $\epsilon>0$ being sufficiently small, that there exist a nonnegative function $e(x)\in H^1(\Omega)$ and a positive constant $t_0$ such that $J_{\epsilon,\delta}(t_0e)=0$.  To this end, we choose a test function following (2.10) in \cite{LNT}.

Without loss of generality, we assume that $\Omega$ contains the origin (otherwise, we can shift $\Omega$ without changing the shape of this test function) and define
\begin{equation}\label{20}
e(x)=
\left\{
\begin{array}{ll}
\epsilon ^{-N}  \left(1-\frac{\vert x \vert}{\epsilon}\right), &\vert x \vert \leq \epsilon, \\
0, &\lvert x \rvert > \epsilon.
\end{array}
\right.
\end{equation}
First of all, we have the following lemma which estimates the energy of this test function.
\begin{lemma}\label{lem21}
There exist two positive constants $C_1$ and $C_2$ with $C_1<C_2$ which are independent of $\epsilon$ and $\delta$ such that
$\frac{dJ_{\epsilon,\delta} (te)}{dt}<0$ if $t>C_1\epsilon^N$ and $J_{\epsilon,\delta} (te)<0$ if $t>C_2\epsilon^N$.
Moreover
\[0<\sup_{t\in (0,\infty)} J_{\epsilon,\delta}(te) \leq C_2 \epsilon^N. \]
\end{lemma}
\begin{proof}We apply the arguments in Lemma 2.4 in \cite{LNT}.  To this end, we shall need to verify that $f_\delta(t)$ in (\ref{21}) satisfies (\emph{h3}) in \cite{LNT}:
\begin{equation}\label{16}
f_\delta(t)/t \rightarrow \infty \text{~as~} t \rightarrow \infty,
\end{equation}
and there exists some positive constants $a_1$ and $a_2$ such that
\begin{equation}\label{19}
f_\delta(t) \leq a_1+a_2 t^p,t\geq0,
\end{equation}
with $p\in(1,\infty)$ if $N=1,2$ and $p\in(1,(N+2)/(N-2))$ if $N\geq3$.  In particular, we need to show that $a_1$ and $a_2$ are uniform in $\delta$.

We can easily see that (\ref{16}) holds for $f_\delta(t)$ uniformly in $\delta$.  To show (\ref{19}), we put $t_\delta=\frac{t_1+c}{\delta^\frac{1}{p-1}}$ and write $f_\delta(t)$ in (\ref{21}) as
\[f_\delta(t)=(t+t_\delta)^p-pt_\delta^{p-1}t-t_\delta^p,\]
 where we put $m=1$ without loss of our generality.  We see that $f_\delta(t)/t^p \rightarrow 1$ for all $\delta\geq \delta_0$ as $t\rightarrow \infty$, therefore, there exists $t_0>0$ independent of $\delta$ such that $f_\delta(t)<2t^p$ for all $t\in(t_0,\infty)$.  Now we take $a_1=\max_{t\in(0,t_0)} f_\delta(t)$, then $a_1$ is uniformly bounded for all $\delta\geq \delta_0$ and we have that
\[f_\delta(t)\leq a_1+2t^p,t\geq0.\]
Therefore this verifies (\emph{h3}) and Lemma \ref{lem21} follows from Lemma 2.4 in \cite{LNT}.
\end{proof}

\begin{lemma}\label{lem22}
There exists a positive constant $\theta \in \big(0,\frac{1}{2}\big)$ independent of $\delta$ such that
\begin{equation}\label{23}
F_\delta(t)=\int_0^t f_\delta(s)ds \leq \theta f_\delta(t)t, \forall t \in(0,\infty).
\end{equation}
Moreover, for any solution $w_{\epsilon,\delta}(x)$ of (\ref{15}), we have that
\begin{equation}\label{24}
\left(\frac{1}{2}-\theta\right)\lVert w_{\epsilon,\delta} \rVert^2_\epsilon \leq J_{\epsilon,\delta}(w_{\epsilon,\delta}) \leq \frac{1}{2}\lVert w_{\epsilon,\delta} \rVert^2_\epsilon.
\end{equation}
\end{lemma}

\begin{proof} We assume $m=1$ as above without loss of our generality.  Denoting $t_\delta=\frac{t_1+c}{\delta^\frac{1}{p-1}}$, we can write $f_\delta(t)$ in (\ref{16}) as
 \[f(t_\delta;t)=(t+t_\delta)^p-pt_\delta^{p-1} t-t_\delta^p,\]
then its antiderivative reads
\[F(t_\delta;t)=\int_0^t f(t_\delta;s)ds=\frac{1}{p+1}(t+t_\delta)^{p+1}-\frac{p}{2}t_\delta^{p-1} t^2-t_\delta^p t-\frac{1}{p+1}t_\delta^{p+1}.\]
Since $\delta>\delta_0$ in (\ref{12}) and $t_1$ is monotone decreasing in $\delta$, therefore $t_\delta\in [0,t^*]$ for all $\delta \in [\delta_0,\infty]$, where $t^*=\frac{t_1+c}{\delta_0^\frac{1}{p-1}}$.

In order to show (\ref{23}) in particular the independence of $\theta$ on $\delta$, it is sufficient to prove that \[\sup_{t_\delta\in[0,t^*]}\max_{t\in(0,\infty)} F(t_\delta;t)-\frac{1}{2}f(t_\delta;t)t=\sup_{t_\delta\in[0,t^*]}\max_{t\in(0,\infty)} g_\delta(t)<0.\]
We first claim that for any $t_\delta\in[0,t^*]$, $g_\delta(t)=F_\delta(t_\delta;t)-\frac{1}{2}f_\delta(t_\delta;t)t<0$ for all $t\in(0,\infty)$.  To show this, we note that $g(0)=0$ and $g'(t)=\frac{1}{2}((t+t_\delta)^p-p(t+t_\delta)^{p-1}t-t_\delta^p )$; moreover, $g''(t)=-\frac{p(p-1)}{2}(t+t_\delta)^{p-2}t<0$ for all $t>0$, therefore $g'(t)$ is monotone decreasing in $t$ and this, together with the fact $g'(0)=0$, implies that $g'(t)<0$ for all $t\in(0,\infty)$.  Therefore, $g(t)$ is decreasing in $t$ and $g(t)<g(0)=0$ for all $t\in(0,\infty)$ as claimed.  Moreover, if $t_\delta=0$ or $\delta=\infty$, we easily see that $F(t_\delta;t)=\frac{1}{p+1}t^{p+1}=\frac{1}{p+1}f(t_\delta;t)t<\frac{1}{2}f(t)t$ for all $t\in(0,\infty)$ since $p>1$.

It is easy to see that $\max_{t\in(0,\infty)}F(t_\delta;t)-\frac{1}{2}f(t_\delta;t)t$ is continuous in $\delta$, therefore $\sup_{t_\delta\in[0,t^*]}\max_{t\in(0,\infty)}F(t_\delta;t)-\frac{1}{2}f(t_\delta;t)t<0$ as desired.  Thus $\theta$ is independent of $\delta$ and this verifies inequality (\ref{23}).

In order to prove (\ref{24}), we first see that the second inequality holds thanks to (\ref{17}) and the definition of $F_\delta(t)$.  To show the first inequality, we test (\ref{15}) by $w_{\epsilon,\delta}$ and then integrate it over $\Omega$ by parts
\begin{equation}\label{27}
\int_\Omega \epsilon^2 \vert \nabla w_{\epsilon,\delta} \vert^2 +c_\delta w_{\epsilon,\delta}^2 dx =\int_\Omega f_\delta(w_{\epsilon,\delta})w_{\epsilon,\delta} dx.
\end{equation}
In light of the energy functional defined in (\ref{17}), we see that
\begin{equation}\label{28}
\begin{split}
J_{\epsilon,\delta}(w_\epsilon)
&\left. =\frac{1}{2}\int_\Omega \epsilon^2 \lvert \nabla w_{\epsilon,\delta} \rvert^2+c_\delta w_{\epsilon,\delta}^2 dx -\int_\Omega F_\delta(w_{\epsilon,\delta}) dx   \right. \\
&\left. \geq \frac{1}{2}\int_\Omega \epsilon^2 \lvert \nabla w_{\epsilon,\delta} \rvert^2+c_\delta w_{\epsilon,\delta}^2 dx-\theta \int_\Omega f_\delta(w_{\epsilon,\delta})w_{\epsilon,\delta} dx,     \right.
\end{split}
\end{equation}
therefore we obtain from (\ref{27}) that
\[J_{\epsilon,\delta}(w_{\epsilon,\delta}) \geq \left(\frac{1}{2}-\theta\right)\int_\Omega \epsilon^2 \lvert \nabla w_{\epsilon,\delta} \rvert^2+c_\delta w_{\epsilon,\delta}^2 dx=\left(\frac{1}{2}-\theta \right) \Vert w_{\epsilon,\delta} \Vert ^2_\epsilon, \]
and this concludes the proof of Lemma \ref{lem22}.
\end{proof}
Now we are ready to present the following existence and nonexistence of nonconstant positive solutions to (\ref{15}).
\begin{proposition}\label{prop2}
Suppose that $\delta>\delta_0$ and $p\in(1,\infty)$ for $N=1,2$ and $p\in(1,(N+2)/(N-2))$ for $N\geq 3$.  Then for each $\epsilon>0$, there exists a positive solution $w_{\epsilon,\delta}(x)$ to (\ref{15}) and the energy functional $J_{\epsilon,\delta}(w_{\epsilon,\delta})$ obtained by
\[J_{\epsilon,\delta}(w_{\epsilon,\delta})=\mathcal{C}_{\epsilon,\delta}=\inf_{v\in H^1(\Omega),v>0} \sup_{t\geq0} J_{\epsilon,\delta}(tv) \]
satisfies
\begin{equation}\label{29}
0<J_{\epsilon,\delta}(w_{\epsilon,\delta})= \mathcal{C}_{\epsilon,\delta}\leq C_2  \epsilon^N,
\end{equation}
where $C_1 $ and $C_2 $ are two positive constants independent of $\epsilon$ and $\delta$.  Moreover, a nonconstant positive $H^1$ solution $w_{\epsilon,\delta}(x)$ to (\ref{15}) is a classical solution and
\begin{equation}\label{30}
0<\inf_\Omega w_{\epsilon,\delta}(x),\sup_\Omega w_{\epsilon,\delta}(x)<C_3
\end{equation}
for some positive constant $C_3$ which is independent of $\delta$.  Furthermore, there exist two positive constants $\epsilon_0$ and $\epsilon^*_0$ which are also independent of $\delta$ such that

(i).  (\ref{15}) has only nonconstant least--energy solution $w_{\epsilon,\delta}(x)$ if $\epsilon<\epsilon_0$, and

(ii).  (\ref{15}) has only constant positive solution $\bar{w} \equiv (t_2-t_1)\delta^{-\frac{1}{p-1}}$ if $\epsilon>\epsilon^*_0$.
\end{proposition}
We shall refer to $w_{\epsilon,\delta}(x)=J^{-1}_{\epsilon,\delta}(\mathcal{C}_{\epsilon,\delta})$ in $H^1(\Omega)$ as a \emph{least--energy solution} of (\ref{15}).
\begin{proof} Since $f_\delta(t)$ satisfies (\ref{16}) and $p$ is subcritical, similar as the proof of Theorem 2 in \cite{LNT}, we can show that the least--energy value $\mathcal{C}_{\epsilon,\delta}$ is achieved at $w_\epsilon(x)$, which is a $H^1$ solution to (\ref{15}).  Moreover, it follows from the standard elliptic regularity argument that $w_\epsilon(x)$ is a classical solution.  On the other hand, similar as the analysis of Section 2 in \cite{LNT}, we can show that $\mathcal{C}_{\epsilon,\delta}\geq C_1\epsilon^N $ and $w_\epsilon(x)$ is uniformly bounded in both $\epsilon$ and $\delta$.  Moreover, nonconstant solution $w_\epsilon(x)$ must be strictly positive over $\bar \Omega$ according to the strong maximum principle and Hopf's boundary lemma.  Furthermore, $\mathcal{C}_{\epsilon,\delta}\leq C_2\epsilon^N $ follows from Lemma \ref{lem21} and the fact that $\mathcal{C}_{\epsilon,\delta}  \leq \sup_{t\geq0} J_{\epsilon,\delta}(te)$ .

The proof of \emph{(ii)} is classical and we will leave it to the reader.  To prove \emph{(i)}, we shall show that $J^{-1}_{\epsilon,\delta}(\mathcal{C}_{\epsilon,\delta})$, which may consist of constant solutions of (\ref{15}), admits only nonconstant positive solutions if $\epsilon$ is small independent of $\delta$.  To this end, we first observe that $\bar{w}\equiv (t_2-t_1)\delta^{-\frac{1}{p-1}}$ is the unique positive solution to (\ref{15}).  To rule out it as a least--energy solution of (\ref{15}) for small $\epsilon$, we have from (\ref{24}) that
\[J_{\epsilon,\delta}(\bar{w})\geq \left(\frac{1}{2}-\theta\right) c_\delta (t_2-t_1)^2 \delta^{-\frac{2}{p-1}} \vert \Omega \vert,\]
however, since $t_2-t_1 >C\delta^\frac{1}{p-1}$ for some $C$ independent of $\epsilon$, we conclude that
\[J_{\epsilon,\delta}(\bar{w})\geq C \vert \Omega \vert,\]
and if $\epsilon$ is small
\[J_{\epsilon,\delta}(\bar{w}) >J_\epsilon(w_{\epsilon,\delta})=C \epsilon^N.\]  Therefore $\bar{w}=t_2-t_1$ can not be a least--energy solution and this finishes the proof of Proposition \ref{prop2}.
\end{proof}

\begin{corollary}\label{cor21}
Let $w_\epsilon(x)$ be a solution of (\ref{15}), then we have that
\[\int_\Omega \epsilon^2 \lvert \nabla w_{\epsilon,\delta} \rvert^2+c_\delta w_{\epsilon,\delta}^2 dx =\int_\Omega f_\delta(w_{\epsilon,\delta}) w_{\epsilon,\delta} dx \leq C   \epsilon^N,\]
\end{corollary}
where $C$ is a positive constant independent of $\epsilon$ and $\delta$.

\begin{proof}
This follows from (\ref{27}) and (\ref{29}).
\end{proof}

\section{Single boundary spike on a platform}\label{section3}

In this section, we construct a positive solution $v_\epsilon(x)$ of (\ref{6}) that has single boundary spike supported on a platform for $\epsilon$ small.  To this end, we first introduce the equation in the entire space, which we shall use to approximate least--energy solution $w_{\epsilon,\delta}(x)$ of (\ref{15}).

\begin{lemma}\label{prop31}
For each $\delta \in (\delta_0,\infty)$, there exists a unique solution $ \mathbf{w}_\delta$ to the following problem
\begin{equation}\label{31}
\left\{
\begin{array}{ll}
\Delta \mathbf{w}-c_\delta \mathbf{w} +f_\delta(\mathbf{w})=0,&x\in \mathbb{R}^N,\\
\mathbf{w}(0)=\max_{x \in \mathbb{R}^N}\mathbf{w}(x)>0,~\mathbf{w}(x)>0,& x\in \mathbb{R}^N,
\end{array}
\right.
\end{equation}
where $c_\delta$ and $f_\delta$ are the same as in (\ref{15}).  Moreover, $\mathbf{w}_\delta$ satisfies the followings:

(i) The solution $\mathbf{w}_\delta$ is radially symmetric such that, $\mathbf{w}_\delta(x)=\mathbf{w}_\delta(\vert x \vert)$;

(ii) $\mathbf{w}'_\delta<0$ for $r>0$, with $r=\vert x \vert$;

(iii) $\mathbf{w}_\delta(r) \leq C e^{-\mu r}$, $r>1$, where $C$ and $\mu$ are positive constants independent of $\delta$.

Furthermore, the system is associated with an energy functional $I_\delta(\mathbf{w})$,
\begin{equation}\label{32}
I_\delta(\mathbf{w})=\frac{1}{2}\int_{\mathbb{R}^N}\vert \nabla \mathbf{w} \vert^2+c_\delta\mathbf{w}^2 dz-\int_{\mathbb{R}^N} F_\delta(\mathbf{w}) dz,
\end{equation}
which is positive and uniformly bounded in $\delta$.
\end{lemma}

\begin{proof} First of all, we see that for all $\delta\in(\delta_0,\infty)$
\[\left( \frac{f_\delta(t)}{t} \right)' \geq 0, \text{ for all }t\geq 0,\]
then by the celebrated Gidas--Ni--Nirenberg \cite{GNN} symmetry theorem and the uniqueness result from Kwong and Zhang \cite{KZ}, there exists a unique solution $\mathbf{w}_\delta(x)$ to (\ref{31}) which is radially symmetric.  Choosing $\epsilon=1$ in the test function $e(x)$ defined in (\ref{20}), we can show that $I_\delta(\mathbf{w}_\delta)$ is uniformly bounded in $\delta$ by the similar proof as for (\ref{29}), moreover there exists a positive constant $C$ independent of $\delta$ such that
\[\lVert \mathbf{w}_\delta \rVert_{H^1(\mathbb{R}^N)} <C ,\]
Then we can apply the radial lemma of Strauss \cite{S} to obtain that
\[\lvert \mathbf{w}_\delta  \rvert \leq Cr^\frac{1-N}{2} \lVert \mathbf{w}_\delta \rVert_{H^1(\mathbb{R}^N)}, r \geq 1,  \]
where $C$ is a positive constant which is independent of $\delta$.  After applying Proposition 4.1 in \cite{GNN} to (\ref{31}), we have that
\[\lvert \mathbf{w}_\delta  \rvert \leq Cr^\frac{1-N}{2} e^{-\mu r}, \]
and this concludes the proof of the lemma.
\end{proof}

This unique solution $\mathbf{w}_\delta$ is called the \emph{ground state} of (\ref{31}).  To study the least--energy solution $w_{\epsilon,\delta}(x)$ of (\ref{15}), we note that $f_\delta \in C^1(\mathbb{R}^+;\mathbb{R}^+)$ and it satisfies $f_\delta(0)=0$, $f_\delta(t)=o(t)$ as $t \rightarrow0^+$ and $f_\delta(t)=O(t^p)$ as $t \rightarrow +\infty$, where $p\in(1,\infty)$ if $N=1,2$ and $p\in(1,\frac{N+2}{(N-2)})$ if $N\geq3$.  By the same arguments in the proof of Theorem 2.3 in \cite{NT} and Theorem 1.2 in \cite{NT2}, we can show the following results on the limiting profiles and single boundary spike of $w_{\epsilon,\delta}(x)$.
\begin{proposition}\label{thm31}%(\cite{NT,NT2})
Let $w_{\epsilon,\delta}(x)$ be a least--energy solution to (\ref{15}), i.e., a critical point of $J_{\epsilon,\delta}$ in $H^1(\Omega)$ such that $J_{\epsilon,\delta}(w_{\epsilon,\delta})=\mathcal{C}_{\epsilon,\delta}$, where $\mathcal{C}_{\epsilon,\delta}$ is the least--energy value.  Then $w_{\epsilon,\delta}(x)$ has at most one local maximum at a point $P_\epsilon$ in $\bar{\Omega}$ if $\epsilon$ is small and $P_\epsilon$ must lie on the boundary $\partial \Omega$ if $\epsilon$ is sufficiently small independent of $\delta$.  For any $\eta>0$, there exist $\epsilon^*(\eta)$ independent of $\delta$ and a small sub--domain $\Omega^{(i)}_\epsilon\subset \Omega$ such that for all $\epsilon \in (0,\epsilon^*)$, the followings are satisfied:

(1). $P_\epsilon \in \partial \Omega^{(i)}_\epsilon\cap \partial \Omega$ and $\text{diam}( \Omega^{(i)}_\epsilon) \leq C \epsilon$;

(2). $\vert w_{\epsilon,\delta} \vert > \eta$ for all $x\in \Omega^{(i)}_\epsilon$, and $\vert w_{\epsilon,\delta} \vert \leq \eta$ for all $x\in \Omega \backslash \Omega^{(i)}_\epsilon$;

(3). $w_{\epsilon,\delta}(x)=\mathbf{w}_\delta\big(\frac{\vert x-P_\epsilon \vert}{\epsilon}\big)+o(\epsilon)$, where $\mathbf{w}_\delta$ is the unique ground state of (\ref{31}).

(4).  Let $H(P)$ be the mean curvature of $P$ at $\partial \Omega$, then
\[H(P_\epsilon) \rightarrow \max_{\bar{\Omega}}  H(P), \text{ as }\epsilon \rightarrow 0.\]
\end{proposition}

Now we proceed to find a positive constant $\delta=\delta_\epsilon$ such that the positive solution $v_{\epsilon,\delta}(x)$ of (\ref{10}) satisfies the integral constraint
\begin{equation}\label{33}
\delta_\epsilon=\int_\Omega (v_{\epsilon,\delta_\epsilon}+c)^p dx.
\end{equation}
To prove (\ref{33}), we first see that it is equivalent to find a $\delta_\epsilon$ such that
\[ \int_\Omega v_{\epsilon,\delta_\epsilon} (x) dx =m.\]
Indeed, we integrate the $v$--equation in (\ref{10}) over $\Omega$ and collect that
\[ 0=\int_\Omega \epsilon^2 \Delta v_{\epsilon,\delta_\epsilon} -v_{\epsilon,\delta_\epsilon} +\frac{m (v_{\epsilon,\delta_\epsilon}+c)^p}{\delta} dx =-\int_\Omega v_{\epsilon,\delta_\epsilon}(x) dx+ \frac{m} {\delta} \int_\Omega(v_{\epsilon,\delta_\epsilon}+c)^p dx, \]
 which shows that (\ref{33}) is equivalent to $\int_\Omega v_{\epsilon,\delta_\epsilon} dx=m$ as claimed.  We now present the proof of our main result.

\begin{proof} [Proof\nopunct] \emph{of Theorem} \ref{theorem11}.  Let $w_{\epsilon,\delta}(x)$ be a least--energy solution obtained in Proposition \ref{thm31}.  Following \cite{SWW}, we denote the set of least--energy solutions to (\ref{15}) as
\begin{equation}\label{34}
\mathcal{S}_\delta=\lbrace  w_{\epsilon,\delta}~\vert~ w_{\epsilon,\delta} \text{ solves } (\ref{15}), J_{\epsilon,\delta}(w_{\epsilon,\delta})=\mathcal{C}_{\epsilon,\delta}  \rbrace.
\end{equation}
It is easy to see that $\mathcal{S}_\delta$ is nonempty; moreover $\mathcal{S}_\delta$ is compact since
\[\sup_{w_{\epsilon,\delta} \in \mathcal{S}_\delta } \lVert w_{\epsilon,\delta} \rVert_\epsilon \leq C,\]
where $C$ is uniform in $\epsilon$ and $\delta$.  Taking
\[\rho(\delta)=\inf_{w_{\epsilon,\delta} \in \mathcal{S}_\delta} \int_\Omega v_{\epsilon,\delta}(x) dx, \]
we have from (\ref{14}) that
\[\rho(\delta)=\delta^\frac{1}{p-1} \inf_{w_{\epsilon,\delta} \in \mathcal{S}_\delta}  \int_\Omega w_{\epsilon,\delta} (x) dx + t_1 \cdot \lvert \Omega \rvert,\]
where $t_1$ is given in Proposition \ref{prop2}.  By the compactness of $\mathcal{S}_\delta$, $\inf_{ w_{\epsilon,\delta} \in \mathcal{S}_\delta}  \int_\Omega w_{\epsilon,\delta}(x) dx$ is a well--defined and continuous function of $\delta$, hence $\rho(\delta)$ is a continuous function of $\delta$.  Furthermore we have from Proposition \ref{thm31} that for any $w_{\epsilon,\delta} \in \mathcal{S}_\delta$,
\begin{equation}\label{35}
\int_\Omega w_{\epsilon,\delta}(x)dx= \epsilon^N\left(\frac{1}{2}\int_{\mathbb{R}^N} \mathbf{w}_\delta(z) dz +o(1)\right),
\end{equation}
where $\mathbf{w}_\delta$ is the unique ground state of (\ref{31}) and $o(1)$ is independent of $\delta$. $\int_{\mathbb{R}^N} \mathbf{w}_\delta(x)\break dz$ is uniformly bounded in $\delta$.  According to Proposition \ref{prop2}, we have that
\begin{equation}\label{36}
\rho(\delta_0)= t_1(\delta_0)\cdot \lvert \Omega \rvert=\frac{c\lvert \Omega \rvert}{p-1} > m,
\end{equation}
where the inequality is due to the assumptions on $M$ and $c$.  On the other hand, since $t_1(\delta) \rightarrow 0$ as $\delta \rightarrow \infty$ uniformly in $\epsilon$, we can find $\delta_1$ large such that
\[t_1(\delta_1) \cdot \lvert \Omega \rvert = \frac{m}{2},\]
then it follows from (\ref{14}) and (\ref{35}) that
\begin{equation}\label{37}
\begin{split}\rho(\delta_1)
&\left. =\int_\Omega v_{\epsilon,\delta_1}(x)dx =\delta_1^\frac{1}{p-1}\int_\Omega w_{\epsilon,\delta_1}(x)+t_1(\delta_1)dx \right.\\
&\left. ={\delta_1}^{\frac{1}{p-1}} \epsilon^N \left(\frac{1}{2}\int_{\mathbb{R}^N} \mathbf{w}_{\delta_1} dz +o(1)\right) +\frac{m}{2}.\right.
\end{split}
\end{equation}
Therefore, we can take $\epsilon$ small enough independent of $\delta$ such that
\begin{equation}\label{38}
\rho(\delta_1)<m.
\end{equation}
Together with (\ref{36}) and (\ref{38}), we observe from the Intermediate value theorem that there exists $\delta_\epsilon \in (\delta_0,\delta_1)$ such that
\[\rho(\delta_\epsilon)=\int_\Omega v_{\epsilon,\delta_\epsilon} (x) dx =m, \text{ thus }\delta_\epsilon=\int_\Omega (v_{\epsilon,\delta_\epsilon}(x)+c)^p dx.\]
Now by taking $\delta=\delta_\epsilon$, we see that Proposition \ref{thm31} implies all but part ($ii$) and (\ref{7}) of Theorem \ref{theorem11}, where $v_{\epsilon,\delta_\epsilon}(x)$ is replaced by $v_\epsilon(x)$.

To prove that the maximum of $v_\epsilon(x)$ has a positive lower bound, we have that
\[v_\epsilon(P_\epsilon)=\max_{\bar\Omega} v_\epsilon(x)>\bar{v}_\epsilon=\frac{1}{\vert \Omega \vert}\int_\Omega v_\epsilon (x) dx=\frac{\beta M}{ \alpha \vert \Omega \vert},\]
where the last identity follows from the integral constraint
\[\int_\Omega v_\epsilon (x) dx= m=\frac{\beta M}{\alpha}.\]

To show (\ref{7}), we put $t^*_\epsilon=t_1(\delta_\epsilon)$ and conclude that $t^*_\epsilon=t_1(\delta_\epsilon)<t_1(\delta_0)=\frac{c}{p-1}$, where the inequality and the last identity follows from Proposition \ref{prop2} and its remark.  On the other hand, we integrate $v_\epsilon(x)$ over $\Omega$ and collect that
\begin{equation}\label{39}
\begin{split}
\int_\Omega v_\epsilon(x) dx
&\left. =\frac{\beta M}{\alpha} = \delta_\epsilon^\frac{1}{p-1}\int_\Omega w_{\epsilon,\delta_\epsilon}(x)+t_\epsilon^*dx \right. \\
&\left. =\epsilon^N \delta_\epsilon^\frac{1}{p-1}\Big(\frac{1}{2}\int_{\mathbb{R}^N} \textbf{w}_{\delta_\epsilon}(z)dz+o(1)\Big)+t_\epsilon^* \cdot \vert \Omega \vert,     \right.
\end{split}
\end{equation}
where the last identity follows from (\ref{35}) or Proposition \ref{thm31} with $\delta=\delta_\epsilon$.  Sending $\epsilon$ to zero in (\ref{39}), we can readily conclude that $t_\epsilon^*$ approaches to $\frac{\beta M}{\alpha \vert \Omega \vert}$ as claimed.  Thus we have proved (\ref{7}) and this completes the proof of Theorem \ref{theorem11}.
\end{proof}

\section{Numerical simulations of spiky solutions}\label{section4}

In this section, we present numerical results on the the formation and evolution of boundary spike of (\ref{2}) over $\Omega=(0,1)\times (0,1)$ to illustrate our theoretical results.  Putting $\alpha=\beta=1$, and $\phi(v)=\ln (v+c)$, we perform extensive numerical simulations on (\ref{2}) by choosing different sets of values for $\epsilon$, $c$ and the initial data.  We shall see that multi--spike solutions also arise through the system.  It is worthwhile to mention that rigorous analysis is needed to fully understand the dynamics such as the large time behaviors and stability of these structures, which is beyond the scope of this paper.

\begin{figure}[ht]
\centering
\includegraphics[width=\textwidth]{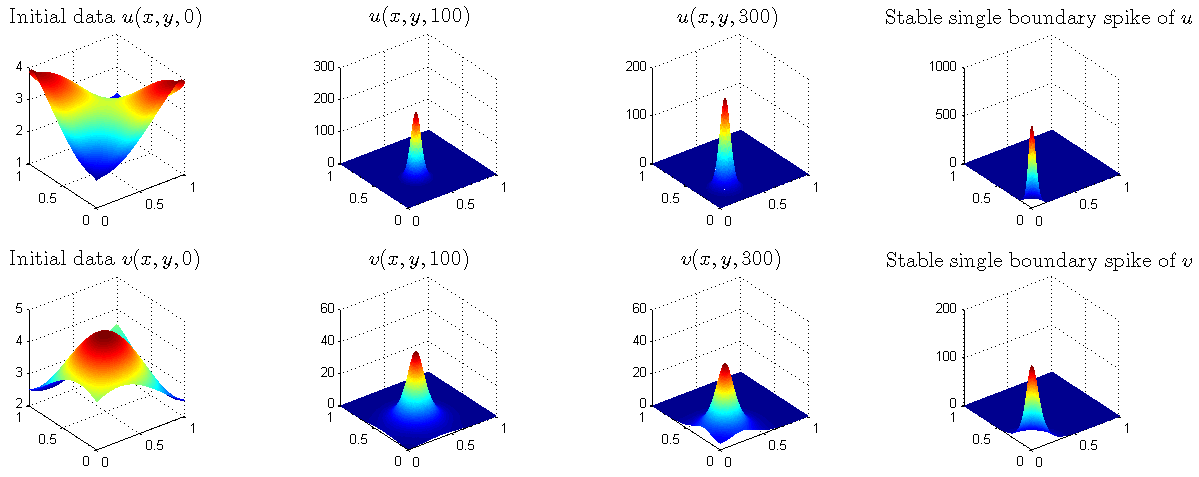}
\caption{The formation of boundary spike of system (\ref{2}) at corner $(0,0)$ of the unit square in 2--D, where initial data are taken to be $u_0(x,y)=3-\cos(\pi x)\cos(\pi y)$, $v_0(x,y)=3+\cos(\pi(x-1/4))\cos(\pi(y-1/4))+\cos(\pi(x-1/2))\cos(\pi(y-1/2))$.  The parameters are $d_1=1$, $d_2=0.01$, $\chi=3$, $\alpha=\beta=1$ and $c=0.1$.}\label{fig1}
\end{figure}

\begin{figure}[ht]
        \centering
                \includegraphics[width=\textwidth]{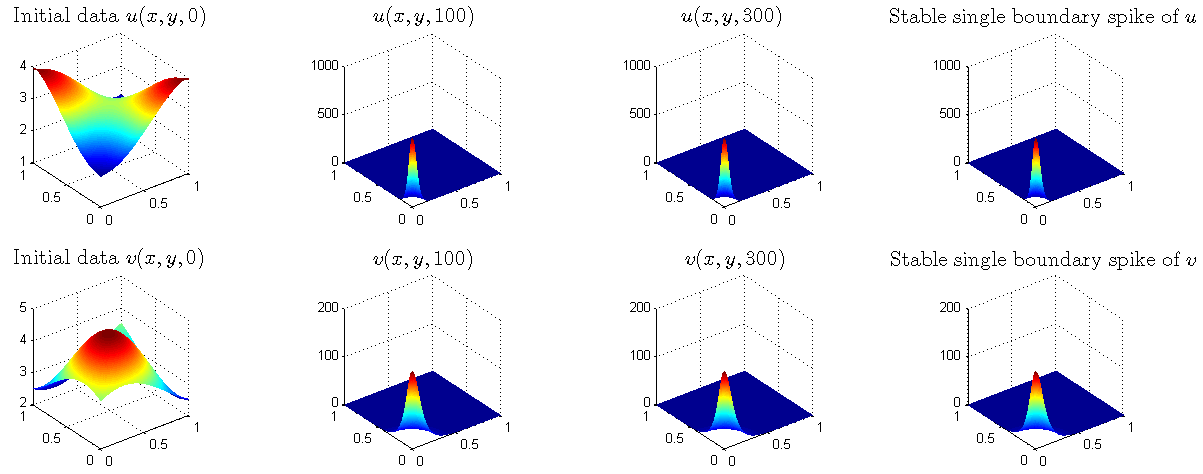}
  \caption{ The initial data $(u_0,v_0)$ and parameters are chosen to be the same as in Figure \ref{fig1} except that $c=5$.  We see that a stable single boundary spike develops at $(0,0)$.
  }\label{fig2}
\end{figure}

\begin{figure}[ht]
        \centering
                \includegraphics[width=\textwidth]{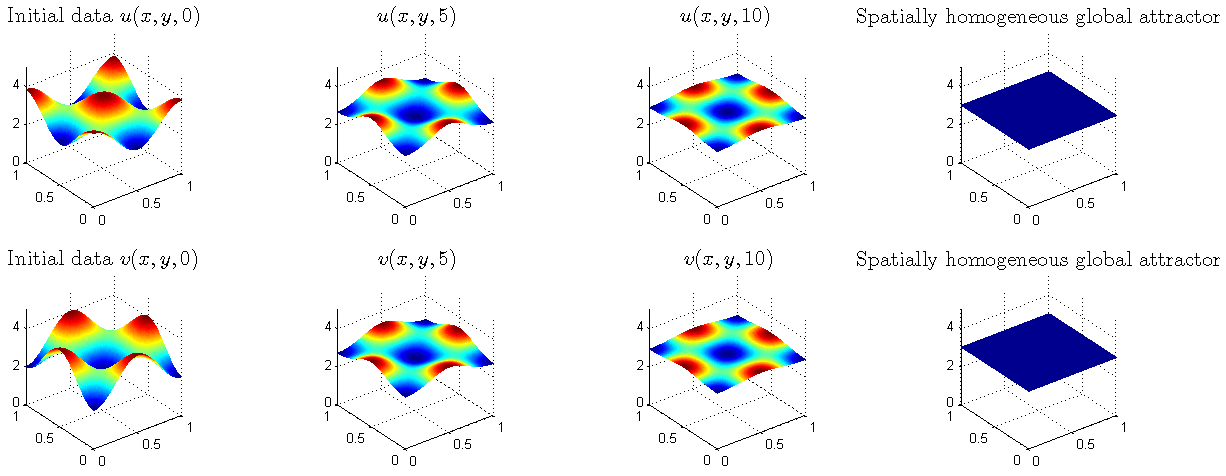}
  \caption{Parameters are chosen to be the same as in Figure \ref{fig1} except that $c=10$.  Initial data are $u_0=3+\cos 2\pi x \cos 2\pi y$ and $v_0=3-\cos 2\pi x \cos 2\pi y$.  We see that solutions converge to the global attractor $(3,3)$.
  }\label{fig3}
\end{figure}

For $c=0.1$, numerical simulations of system (\ref{2}) on the unit square are plotted in Figure \ref{fig1}.  The graphes in the first line represent the spatial--temporal behaviors of the cellular population density $u(x,y,t)$ and the graphes in the second line illustrate the behaviors of the chemical concentration $v(x,y,t)$.  Our simulations show that $u$ quickly develops an interior spike which gradually moves to the corner.  Eventually a stable spike is formed at the boundary point corner $(0,0)$, where the boundary is most curved.
In Figure \ref{fig2}, we choose $c=5$ and take the same initial data and system parameters as in Figure \ref{fig1}.  For this set of values, solutions quickly evolves into a boundary spike at $t=100$ compared to that in Figure \ref{fig1}.  It eventually develops into a stable boundary spike concerned at $(0,0)$.  We find through our numerical results that the stable boundary spike stays on a platform with height 0.02, compared with those obtained in Figure \ref{fig1}.  According to our steady state analysis, we surmise that this platform is contributed by the large magnitude of $c$.  We need to point out that our numerical results are ambiguous since the height of the platform for Figure \ref{fig2} is $t^*=\frac{\beta M}{\alpha \vert \Omega \vert}=3$, which is greatly larger than 0.02.  However, the numerical simulations support our analysis qualitatively and rigorous stability analysis is needed to fully understand the dynamics of the spiky solutions.  We also perform numerical simulations in Figure \ref{fig3} which show that the constant solution is a global attractor of (\ref{2}) if $c$ is sufficiently large, independent of the initial data.  This corresponds to the fact that the saturating coefficient $c$ in (\ref{6}) reduces the chemo--attraction magnitude and stabilizes the positive constant steady state.

\begin{figure}[ht]
\centering
\includegraphics[width=\textwidth]{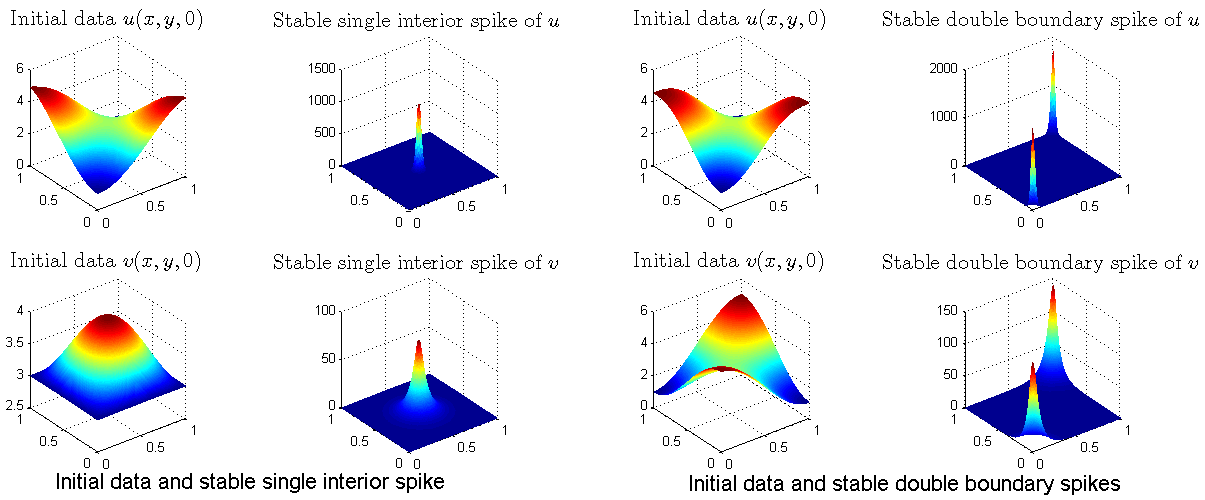}
\caption{Two sets of initial data develop into a stable single interior spike and a stable double boundary spike respectively.  Parameters are chosen to be $d_1=1$, $d_2=0.01$, $\chi=5$, and $\alpha=\beta=c=1$.  Our results indicate that spikes usually develop at the spots where initial chemical concentration maximizes.}\label{fig4}
\end{figure}

We proceed to simulate the formations of stable interior spike and multi--spikes in Figure \ref{fig4}.  Our goal is to illustrate the effect of initial data on the localization of the spikes, i.e., the spot where a spike is formed.  To this end, we choose different initial data in Figure \ref{fig4}.  The initial data on the left are taken to be $u_0(x,y)=3-\cos(\pi x)\cos(\pi y)-\cos(\pi (x-1))\cos(\pi (y-1))$, $v_0(x,y)=3+\cos(\pi(x-1/2))\cos(\pi(y-1/2))$, where the initial chemical concentration $v_0$ attains its maximum at point $(1/2,1/2)$.  We see that a stable interior spike is formed at this point.  The initial data on the right are $u_0(x,y)=3-\cos(\pi x)\cos(\pi y)-\cos(\pi (x-1))\cos(\pi (y-1))$, $v_0(x,y)=3+\cos(\pi x)\cos(\pi y)+\cos(\pi (x-1))\cos(\pi (y-1))$ and $v_0$ achieves its maximum at two boundary points $(0,0)$ and $(1,1)$.  Then a stable double boundary spike is developed at both maximum points.  Our numerical simulations indicate that the spike of (\ref{2}) with saturated logarithmic sensitivity can localize at the spot where the initial chemical concentration reaches its largest value, independent of the initial cellular population distribution.  This complies with the well--accepted theoretical results that chemotaxis dominates the dynamics of the cellular movements over the domain when $\chi$ is not too small or $c$ is not too large.  However, a rigorous analysis for this purpose requires totally non--trivial mathematical techniques.

\begin{figure}[ht]
        \centering
\includegraphics[width=\textwidth]{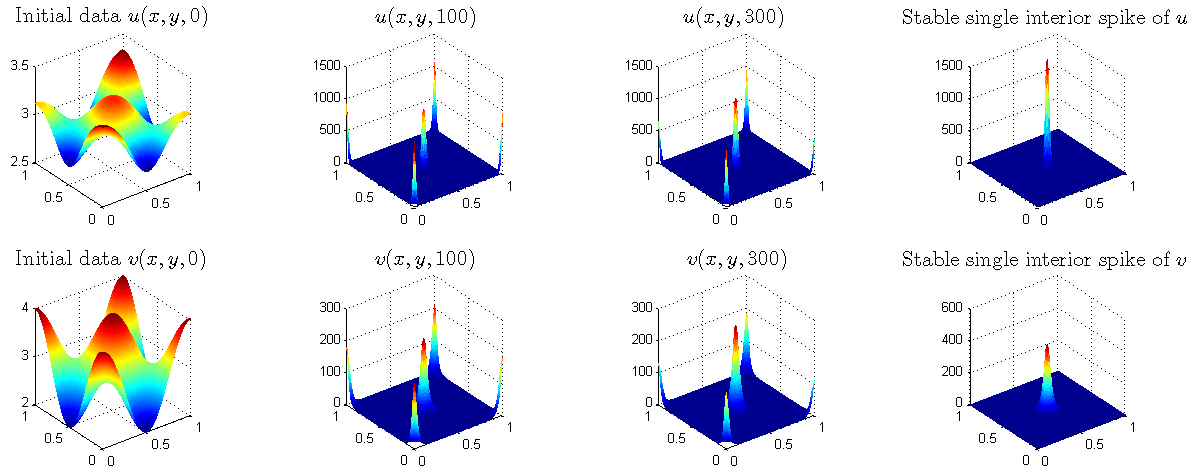}
        \caption{The formation and evolution of multi--spikes.  The initial data are $u_0(x,y)=3+0.1\cos(\pi x)\cos(\pi y)$, $v_0(x,y)=3+\cos(2\pi x)\cos(2\pi y)$.  The parameters are chosen $d_1=1$, $d_2=0.001$, $\chi=3$ and $\alpha=\beta=c=1$.}\label{fig5}
\end{figure}

Finally, we present Figure \ref{fig5} to illustrate the emergence of multi--spikes in model (\ref{2}) with parameters taken to be $d_1=1$, $d_2=0.001$, $\chi=3$ and $\alpha=\beta=c=1$.  The initial data are $u_0(x,y)=3+0.1\cos(\pi x)\cos(\pi y)$ and $v_0(x,y)=3+\cos(2\pi x)\cos(2\pi y)$.  We observe that both $u$ and $v$ develop multi--spikes ($t=100$) at the four corners and the center very quickly.  The structures of these multi--spikes keep well--preserved for quite a long time period (through $t=100$ to $t=300$), then four spikes at the corner disappear and both $u$ and $v$ develop into a stable interior spike at the center eventually.  Our numerical simulations illustrate that the multi--spikes can arise from system (\ref{2}) and they are always unstable or meta--stable.  The rigorous stability analysis is beyond the scope of this paper.  Figure \ref{fig5} also suggests that the formation of stable patterns, i.e., solutions with spikes or stripes, etc., depends on the initial data and in particular the initial distribution of chemical concentration.

\section*{Acknowledgments}
This work was initiated when I was a PhD student at Tulane University.  I wish to express my sincere gratitude to my thesis advisor Professor Xuefeng Wang for his guidance and support over the years, as well as his helpful discussions and comments on an earlier version of this paper.  I want to thank the anonymous referee for carefully reading this paper and giving insightful comments and constructive suggestions which improved its exposition.  This work was partially supported by the Summer Research Fund from the Math Department of Tulane University.  This research also receives support from the Project--sponsored by SRF for ROCS, SEM.

\medskip
% The data information below will be filled by AIMS editorial staff
%Received May 2014; revised  November 2014.
\medskip
\end{document}